**UDK 51(091)**

# On the history of epsilontics

### G.I. Sinkevich

St. Petersburg University of Architecture and Civil Engineering

4 Ul. 2-ya Krasnoarmeiskaya, St. Petersburg

*e-mail*: galina.sinkevich@gmail.com

This is a overview of the genesis of ε–δ language in works of mathematicians of the 19th century. It shows that although the symbols ε and δ were initially introduced in 1823 by Cauchy, no functional relationship for δ as a function of ε was ever specified by Cauchy. It was only in 1861 that the epsilon-delta method manifested itself to the full in Weierstrass' definition of a limit. The article gives various interpretations of these issues later provided by mathematicians.

> It is mere feedback-style ahistory to read Cauchy
> (and contemporaries such as Bernard Bolzano)
> as if they had read Weierstrass already.
> On the contrary, their own pre-Weierstrassian
> muddles need historical reconstruction.
> *I. Grattan-Guinness* [Grattan Guinness, 2004, p. 176].

Since the early ancient word, the concept of continuity was described throw notions of time, motion, divisibility, contact[1].

With the extension of the range of tasks and the development of ideas about the notion of a function, the understanding of continuity on the basis of physics and geometry became insufficient, ultimately necessitating an arithmetization of this concept.

In the 17th century, G. Leibnitz laid down The Law of Continuity as follows: "In any supposed [continuous] transition, ending in any terminus, it is permissible to institute a general reasoning, in which the final terminus may also be included[2]" [Child, p.40].

---

[1] "The 'continuous' is a subdivision of the contiguous: things are called continuous when the touching limits of each become one and the same and are, as the word implies, contained in each other: continuity is impossible if these extremities are two". "Being continuous it is one" [Aristotle, Physics, Book V, parts 3, 4].
[2] The word "continuous" was omitted is Child's translation, as pointed

In his *Arithmetica Infinitorum* Wallis introduced the following definition of a limit: "But (which for us here suffices) they continually approach more closely to the required ratio, in such a way that at length the difference becomes less than any assignable quantity" [Wallis, 1656, p.42]. Euler's personal copy of this book of Wallis' is preserved in the Archives of St. Petersburg Academy of Sciences as part of the Euler collection.

Euler considered functions represented by a single formula to be continuous. According to Euler, "the rules of calculation are based on The Law of Continuity, pursuant to which curved lines are drawn by continuous movement of a point", "a continuous line is drawn so that its nature is presented with the help of one specific function of *x*» [Euler, L. 1748 (1961), v.2, p. 21].

In 1765, J. d'Alembert provided the following definition of a limit: "A value is said to be a limit of another value if the latter can approximate the former nearer than any given value, no matter how small it may be supposed, however, without the approximating value being able to exceed the value it approximates; thus, the difference between such value and its limit is absolutely indeterminable" [D'Alembert, 1765, p. 155-156]. The limit in the d'Alembert was not constant.

L. Carnot in 1797 tried to unite the method of limits and infinitesimal calculus for "procedures of both methods became absolutely identical" [Yushkevich, 1986, p. 45-55]. The contest declared by Berlin Academy of Sciences in 1786 at the initiative of G. Lagrange promoted strengthening of interest to infinitesimal issues: "... we need a *clear and accurate theory of what is called continuous in mathematics"* [Yushkevich, 1973, p. 140]. None of the 23 works submitted to the contest satisfied the Academy: "…the principle we need must not be limited to calculation of infinitely small values; it must extend to algebra and geometry as well which render after the manner of the ancients" [ibib., p. 141]. The winner of the contest was a Swiss mathematician Simon l'Huilier (1750 – 1840) who lived in Warsaw at that time. It was in his work entitled 'Elementary statement of principles of calculus' published by the Academy in 1786 that symbol out by M. Katz and D. Sherry [Katz, M., Sherry D., 2012, p. 1551]. I am grateful to S.S. Demidov who brought to my attention the similarity of this idea of Leibnitz and Bolzano's and Cauchy's understanding of continuity [Demidov, 1990].

$\lim \dfrac{\Delta P}{\Delta x}$ first appeared [L`Huilier, 1786, p. 31]. Thereafter, this symbol was used by Lacroix[3].

The infinitesimal methods disappointed Lagrange and in subsequent years, he avoided using infinitely small values (infinitesimals). Nonetheless, in 1811, in his second edition of Analytical Mechanics, Lagrange called infinitesimals a safe and convenient tool to simplify a proof[4].

"On a conservé la notation ordinaire du calcul differentiel, parce qu'elle répond au système des infiniment petits, adopté dans ce Traité. Lorsqu'on a bien conçu l'esprit de ce systeme, et qu'on s'est convaincu de l'exactitude de ses résultats par la méthode géométrique des premières et dernières raison, ou par employer les infiniment petits comme un instrument sur et commode pour abréger et simplifier, les démonstrations. C'est ainsi qu'on abrège les démonstrations des Anciens, par la méthode des indivisibles » [Lagrange, 1811, p.ii].

The most popular method of 18[th] century geometricians was approximation. For example, "solving an equation like $(x+1)^{\mu} = a$ with noninteger $\mu$, we cannot come up with a precise solution, but approximate it using infinite series. Having determined the finite number of elements of the approximating series, 18[th] century geometricians tried to calculate the upper limit of the approximation error (error, ε), i.e. the difference between the sum of series and its *n-th* subsum. The inequation algebra served as evidence-based techniques here" [Grabiner, 1983, p. 4 of the electronic version].

First decades of the 19[th] century can be characterized as a period of "naïve" theory of functions – analysis developed on the basis of elementary functions, both continuous and differentiable, based on intuitive qualitative definitions of a limit, neighborhood, continuity, and convergence.

In 1797, Lagrange published 'The Theory of Analytical Functions which contains basics of differential calculus free from any consideration of infinitely small limits tending to zero and fluxions, and reduced to analysis of finite values'. Considering function $fx$ and substituting $x$ with a new value $x + i$, Lagrange asserts that $f(x+i)$ can be expanded into series by positive degree $i$, their coefficients being determined by way of differentiation, which is true for known functions. Considering the first expansion term, Lagrange obtains $f(x+i) = fx + iP$, from which $P = \dfrac{f(x+i) - fx}{i}$. In this case, $i$ can be so small that any expansion term would be greater than the sum of all subsequent expansion terms, and the same applies to all smaller $i$ values as well [Yushkevich A. 1977, p. 160–168]. Lagrange adds: "Perfection of approximation methods in which series are applied depends not only on convergence of series, but on the ability to assess the error resulting from neglected terms as well; and it can be stated that all approximative methods used in geometrical and mechanical problems are so far very imperfect. The previous theorem will in many cases be able to inform of the perfection they miss, failing which they are often dangerous to apply" [Lagrange, p. $67 - 68$][5].

---

[3] Sylvestre Lacroix (1765 – 1843) was Lagrange's follower in Ecole Politechnique and professor of analysis for Cauchy. In 1850s, Weierstrass started using symbol $\lim\limits_{x=c}$ ; in 1905, an English mathematician John Lesem first used $\lim\limits_{x \to c}$ in his book [Leathem, 1905].

[4] «On a conservé la notation ordinaire du calcul differentiel, parce qu'elle répond au système des infiniment petits, adopté dans ce Traité. Lorsqu'on a bien conçu l'esprit de ce systeme, et qu'on s'est convaincu de l'exactitude de ses résultats par la méthode géométrique des premières et dernières raison, ou par employer les infiniment petits comme un instrument sur et commode pour abréger et simplifier, les démonstrations. C'est ainsi qu'on abrège les démonstrations des Anciens, par la méthode des indivisibles» [Lagrange, 1811, p.ii].

In 1800, appears C.F. Gauss' work entitled 'Basic concepts of the theory of series' (see [Gauss,1800]) where he considers series as sequences of subsums.

In 1806, André-Marie Ampère published his article entitled 'Elaboration of certain issues in differential calculus which enable obtaining a new presentation of Taylor expansion and expression thereof of closed form if summing is limited' [Ampére, 1806] which is directly relevant to our story. Here, Ampère proves Lagrange mean value theorem on 33 pages and, based thereon, obtains what we know as Taylor expansion with an integral reminder in the form of Lagrange. A.P. Yushkevich calls this Ampère's work an attempt to prove continuous function differentiability analytically [Yushkevich, 1972, p. 243].

Ampère's key tool of proof were inequations[6]. Using them, he assessed approximants and characterized interpolation error. Keeping to Lagrange, Ampère considers $\dfrac{f(x+i) - f(x)}{i}$ as a function of two variables $x$ and $i$ which represents a divided difference of two values of $x$ and $x + i$ of the same variable, this difference not equaling zero or infinity at whatsoever $x$, while at $i = 0$ it changes into $\dfrac{0}{0}$, however, equals neither zero nor infinity. Lagrange called this function *resulting from a derivate*.

---

[5] Quoted from [Yushkevich A. 1972, p. 298] as translated by A.P. Yushkevich.
[6] The same method is used in works of G. Lagrange, J.-B. Fourier (1822), and P.A. Rakhmanov (1803).

Note that symbol *i* here identifies a real number; an imaginary unit was denoted by symbol $\sqrt{-1}$ at that time. Ampère explicitly states that he will consider only functions of a real variable. Naturally, the consideration included by default only "good" functions, i.e. continuous functions and those differentiated at the finite interval[7]. Ampère notes that this function must decrease and increase as *i* changes. Variable *x* changes from $x = a$ to $x = k$, the respective values of function $f(x)$ being denoted through *A* and *K*. Ampère divides the interval from $x = a$ to

$x = k$ by intermediate values *b*, *c*, *d*, *e*, which satisfy (correspond to) the values of function *B*, *C*, *D*, *E*. Thereafter, he builds divided differences like $\dfrac{K-E}{k-e}$ and $\dfrac{E-A}{e-a}$ and proves the correctness of inequations like $\dfrac{E-A}{e-a} < \dfrac{K-A}{k-a} < \dfrac{K-E}{k-e}$. Further, new values are inserted between the old ones, and new inequations are written. As a result, for certain *x*, $f'(x)$ gradually approaches the value of $\dfrac{f(x+i)-f(x)}{i}$. It follows from here that this value is always between two values of a derivative calculated between

$$x \text{ and } x+i.$$

Let us assume that $x+i = z$ and $\dfrac{f(z)-f(x)}{z-x} = p$. Then $f(z) = f(x) + p \cdot (z-x)$. Proceeding with the procedure, Ampère obtains

$$f(z) = f(x) + f'(x) \cdot (z-x) + p' \cdot (z-x)^2,$$

$$f(z) = f(x) + f'(x) \cdot (z-x) + \frac{f''(x)}{2} \cdot (z-x)^2 + \frac{p''}{2}(z-x)^3,$$

$$f(z) = f(x) + f'(x) \cdot (z-x) + \frac{f''(x)}{2} \cdot (z-x)^2 + \frac{f'''(x)}{2 \cdot 3}(z-x)^3 + \frac{p'''}{2 \cdot 3}(z-x)^4,$$

etc.

---

Ampère gives examples of expansion of certain elementary functions. Further, considering $f(x)$ to be primary (primitive) relative to $f'(x)$, he finds the way the sign of a derivative depends on the increase or decrease of the function [Ampére, 1806]. Ampère's proof looks lengthy and unskilful. It was this particular inadequacy that made Augustin Louis Cauchy (1789 − 1857) wish to offer a concise and beautiful construction, which, as we will see later, was used to create the ε-δ language.

Since 1813, Cauchy was teaching at Ecole Politechnique; in 1816, he became member of the Academy. In 1821, he published the Course of Analysis [Cauchy, 1821] (translated into Russian − [Cauchy, 1821, 1864]) he had given at Ecole Royal Politechnique. In this Course Cauchy provides a definition of a continuous function as follows: "Function $f(x)$ given between two known limits of variable $x$ is a continuous function of this variable provided that for all values of variable $x$ taken between these limits the numerical value of difference $f(x+\alpha) - f(x)$ infinitely decreases together with $\alpha$. In other words, function $f(x)$ remains continuous for all $x$ in between two given limits if an infinitely small increase in the variable in between these limits would always imply an infinitely small increase in the function itself. We would also add that function $f(x)$ which is continuous for $x$ will be continuous for voisinage of variable $x$ between the same limits, no matter how close to these limits $x$ is" [Cauchy, 1821, p. 43]. Here he understands the limit as the terminal point of the interval concerned.

Henceforth, each time he referred to a continuous function, Cauchy repeated and used this definition only. An English historian of mathematics J. Gray notes: "Cauchy defined what it is for a function to be integrable, to be continuous, and to be differentiable, using careful, if not altogether unambiguous, limiting arguments."[Gray, p. 62]. Błaszcyk, Katz and Sherry pointed out that Gray is not being accurate when he includes continuity among properties Cauchy allegedly defined using "limiting arguments". Namely, the word "limit" does appear in these Cauchy's definition, but only in the sense of endpoint of the interval, not in any sense related to the modern notion of the limit [Błaszcyk, Katz, Sherry, 2012].

In § 3 of Chapter One of the Course of Analysis, Cauchy considers special values of the function and proves a theorem he is going to need for consideration of equivalency of infinitesimals [Cauchy, 1864, p. 46][8]:

[8] Translated by F. Ewald, V. Grigoriev, A. Ilyin.

"If with increase of variable $x$ difference $f(x+1) - f(x)$ tends to known limit $k$, fraction $\dfrac{f(x)}{x}$ will at the same time tend to the same limit as well.

*Proof.* Let us suppose that number $k$ has a finite value and that $\varepsilon$ is an arbitrary small number. According to the statement, as $x$ increases, difference $f(x+1) - f(x)$ tends to limit $k$; besides, one can always take such great number $h$ that at $x$ equal to or greater than $h$ this difference will be constantly between limits $k - \varepsilon$, $k + \varepsilon$. Having assumed this, let us denote any whole number through $n$, then each quantity will be as follows: $f(h+1) - f(h), f(h-2) - f(h+1), ..., f(h+n) - f(h+n-1)$, and therefore, their arithmetical average, i.e. $\dfrac{f(h+n) - f(h)}{n}$, will be in between limits $k - \varepsilon$, $k + \varepsilon$. Therefore, $\dfrac{f(h+n) - f(h)}{n} = k + \alpha$, where $\alpha$ is the number between limits $-\varepsilon$, $+\varepsilon$.

Let us suppose now that $h + n = x$. Then the previous equation will become

$$\frac{f(x) - f(h)}{x - h} = k + \alpha, \qquad (1)$$

Therefore, $f(x) = f(h) + (x - h) \cdot (k + \alpha)$ and

$$\frac{f(x)}{x} = \frac{f(h)}{x} + \left(1 - \frac{h}{x}\right) \cdot (k + \alpha). \qquad (2)$$

For the value of $x$ to be able to increase indefinitely, it is sufficient to increase number $n$ indefinitely without changing value $h$. Therefore, let us assume that $h$ is constant in equation (2) and $x$ is a variable tending to limit $\infty$; then the numbers of $\dfrac{f(h)}{x}, \dfrac{h}{x}$ contained in the second part will tend to limit zero and the entire second part, to a limit that can be described as $k + \alpha$, where $\alpha$ is constantly confined between $-\varepsilon$ and $+\varepsilon$. Therefore, the limit of relation $\dfrac{f(x)}{x}$ will be the number confined between $k - \varepsilon$ and $k + \varepsilon$.

Whereas this conclusion is true no matter how small is $\varepsilon$, the unknown limit of the function will be number $k$. In other words,

$$\lim \frac{f(x)}{x} = k = \lim \left[ f(x+1) - f(x) \right]."$$

The case where $x$ tends to $\pm\infty$ is considered in the same way [Cauchy, 1864, p. 46 – 49].

As we can see, there is a structure here the development of which led to introduction of ε–δ method. However, ε here is a finitesimal, although arbitrary small, assessment of an error. Cauchy improves Ampère's construction. Two years later, he will improve the rationale from this proof. However, so far, the need to read the course in the customary way without detouring to developments did not enable Cauchy to experiment introducing new methods. Judging from the fact that Cauchy had to explain basics (reduction to common denominator, fundamentals of trigonometry, properties of exponentials) to his students, their basic training was quite limited. Students are known to clamor against studying complex numbers which they believed to be an absolutely useless domain of math.

Cauchy's basic course includes statement of elementary terms of one or more variables, terms of real and imaginary variables (complex variable used to be called imaginary variable at that time), their properties, theory of limits including comparison of infinitesimals, theory of series, Lagrange interpolation formulas.

In 1822, J.-B. Fourier's Analytic Theory of Heat where he used δ-changes was published [Fourier, 1822, p.139].

In 1823, Lecture Notes were published based on the course of lectures in Infinitesimal Calculus [Cauchy, 1823] read by Cauchy at Ecole Politechnique. The course was intended for 40 lectures. These notes were published in Russian under the title "Kratkoje izlozhenije urokov o differenzialnom i integralnom ischislenii" (Differential and Integral Calculus) as translated by V.Y. Bunyakovsky in 1831 [Cauchy, (1823),1831]. The book contained a definition of a limit as follows: "If values attributed to any variable number approximate the value determined so as to finally differ from the latter as small as desired, then these former values are called the limit of all others"[9] [ibid., p. 3] and a definition of a continuous function as follows: "If function $f(x)$ changes along with value $x$ so that for each value of this variable confined in the given limits it has one unambiguous value, then difference $f(x+i) - f(x)$ between the limits of value $x$ will be an infinitely small number; while

---

[9] In connection therewith, the remark of an English historian of mathematics, J. Gray: "Although limits did appear in Cauchy's definitions, however, they meant only the finite point of the definition interval" [Gray, 2008, p. 62] seems inappropriate. As Mikhail Katz notes, "Gray claims that Cauchy defined continuity using "limiting arguments". This is inaccurate. In our paper [Błaszcyk, Katz, Sherry, 2012] we point out the inaccuracy of what he wrote, and add that limits did appear in Cauchy's definition of continuity, but only as the endpoint of the interval of definition" [Personal message].

function $f(x)$ which meets this condition is called a continuous function of variable $x$ between those limits" [ibid., p.11]. And further, in the second lecture:

"If variables are linked with one another so that judging from the value of one variable, values of the rest of variables can be obtained, this means that these different values are expressed with the help of one of them called *independent variable,* and the values presented through it are called *functions* of this variable.

Letter $\Delta$ is often used in calculations to denote a concurrent increase of two variables depending on each other[10]. In this case, variable $y$ will be expressed as a function of variable $x$ in equation

$$y = f(x). \tag{3}$$

Therefore, if variable $y$ is expressed as a function of variable $x$ in equation $y = f(x)$, then $\Delta y$, or an increase in $y$ caused by an increase in $\Delta x$ of variable $x$, will be denoted by formula:

$$y + \Delta y = f(x + \Delta x). \tag{4}$$

<…> It is evident that (1) and (2) are interrelated, therefore

$$\Delta y = f(x + \Delta x) - f(x). \tag{5}$$

Now, let $h$ and $i$ be two different values, the first of which is finitesimal and the second one, infinitely small; and let $\alpha = \dfrac{i}{h}$ be infinitely small value given as a relation of these two values. If finite value $h$ corresponds to $\Delta x$, then value $\Delta y$ given in equation (5) will be the so-called finite difference of function $f(x)$ and will naturally be the finite number.

Should you conversely give $\Delta x$ an infinitesimal value, e.g. $\Delta x = i = \alpha h$, the value of $\Delta y$ will be $f(x + i) - f(x)$ or $f(x + \alpha h) - f(x)$, and will naturally be infinitely small. It is easy to verify that in the context of functions $A^x, \sin x, \cos x,$ to which the following differences correspond:

$$A^{x+i} - A^x = \left( A^i - 1 \right) \cdot A^x,$$

---

[10] There was no such remark in the course of 1821. Here Cauchy points to the existence of a link between increment of function and increment of argument without detailing the connection in their change, which was done by Weierstrass forty years later. Instead, there appears term 'concurrent' (simultané) which is characteristic of the 18th and 19th century. Moreover, the exhaustion method was compared to anthropomorphous time. Newton said that he could calculate the area under the parabola over half a quarter of an hour; he also said (see [Cajory F. 1919, p. 103]): "at the moment the hour expires, no inserted or described figure exists anymore; however, each of them aligns with a curvilinear figure which is the limit they reach". Other mathematicians of the 18th century also defined the limiting process as taking some hours, that is, eventually observable. In this event, this symbol ε denoted a calculation error in Cauchy's works as well.

$$\sin(x+i) - \sin x = 2\sin\frac{i}{2}\cos\left(x+\frac{i}{2}\right),$$
$$\cos(x+i) - \cos x = -2\sin\frac{i}{2}\cos\left(x+\frac{i}{2}\right),$$

each of these differences has multiplier $\left(A^i - 1\right)$ or $\sin\frac{i}{2}$ which along with $i$ infinitely approaches to a limit that equals zero.

Thus, for function $f(x)$ which uniquely possesses finite values for all $x$ contained in between the two given limits difference $f(x + i) - f(x)$ will always be infinitely small between these limits, i.e. $f(x)$ is a continuous function within the limits it changes.

Function $f(x)$ is also said to always be a continuous function of variable $x$ in the vicinity of any particular value of such variable, if this function is continuous between two (even quite close) limits containing this given point" [Cauchy, 1823, p. 17].

On the assumption that any continuous function is differentiable, Cauchy proves the mean value theorem (see [ibid., p. 44 − 45]; [Cauchy, (1823), 1831, p. 36]) as follows:

"*Theorem*. Let function $f(x)$ be continuous between two limits $x = x_0, x = X$. Let us denote the greatest value of its derivative through $A$, $B$ being the smallest value of its derivative between the same limits. Then divided difference $\frac{f(X) - f(x_0)}{X - x_0}$ will inevitably confined between $A$ and $B$.

Let us denote infinitesimal numbers with letters δ, ε, of which let the first one be such number that for numeric values of $i$ that are less than δ, and for any value of $x$ confined between limits $x_0$, $x$ relation $\frac{f(x+i) - f(x)}{i}$ will always be greater than $f'(x) - \varepsilon$ and smaller than $f'(x) + \varepsilon$."[11]

Cauchy mentions that in this proof he keeps to Ampère's memoir quoted above.

Like Ampère, Cauchy inserts new values[12] of $x_1, x_2, \ldots, x_{n-1}$ between $x_0$ and $x$ so that difference $X − x_0$ could be decomposed into positive parts $x_1 − x_0$, $x_2 −$

---

$x_1,\ldots,$ $X - x_{n-1}$, which do not exceed $\delta$. "Fractions $\dfrac{f(x_1)-f(x_0)}{x_1-x_0}, \dfrac{f(x_2)-f(x_1)}{x_2-x_1}, \dfrac{f(X)-f(x_{n-1})}{X-x_{n-1}}$, located between limits: first: $f'(x_0)-\varepsilon, f'(x_0)+\varepsilon$, second: $f'(x_1)-\varepsilon, f'(x_1)+\varepsilon$, will be greater than $A-\varepsilon$, however, no smaller than $B+\varepsilon$. Whereas denominators of the fractions have the same sign, having divided the sum of their numerators by the sum of their denominators, we will obtain an intermediate value fraction, that is to say, a fraction the value whereof lies between the smallest and the greatest of the fractions. However, whereas $\dfrac{f(X)-f(x_0)}{X-x_0}$ is an intermediate value fraction, it is therefore confined between limits

$A - \varepsilon$ and $B + \varepsilon$. And whereas it is true at as

small $\varepsilon$ as we please, $\dfrac{f(X)-f(x_0)}{X-x_0}$ therefore lies between limits $A$ and $B$" [9, p. 36] and [27, p. 44][13].

$$A < f'(x)-\varepsilon < \frac{f(x+i)-f(x)}{i} < f'(x)+\varepsilon < B \quad \text{для } i < \delta.$$

In other words,

Cauchy brilliantly simplified Ampère's proof, having introduced simpler symbols. Ampère states his proof on half the 33 pages, while Cauchy's proof is stated on two pages. Ampère introduces eight auxiliary values, and creates an estimation of relation for each value; instead of averaging, he proves lengthy inequations. Cauchy's proof is elegant and concise.

But Cauchy does not analyze interdependence of $\varepsilon$ and $\delta$ or dependence of $\delta$ on the ensuing difference between neighboring values of the variable. Essentially, $\delta$ is included in a proclaimatory way, irrelatively of the rest of the reasoning.

---

[13] Translated by V.Y. Bunyakovsky.

An Amerian researcher, Judith Grabiner believes [Grabiner] that Cauchy transformed proving technique of inequation algebra into a rigorous approximation error assessment tool.

A Dutch researcher, T. Koetsier believes [Koetsier] that Cauchy arrived at his continuity concept analyzing his proof of the mean-value theorem maybe only for polynomials. It is evident that $x_n$ in his proof are variables which differ from an infinitely small by a constant $a$. Pursuant to Cauchy's definition of continuity, $f(x_n)$ must differ from $f(a)$ by an infinitesimal value. Unlike Grabiner, analyzing Cauchy's proof, Koetsier finds no trace of $\varepsilon - \delta$.

Analyzing Grabiner's hypotheses that Cauchy only assessed the approximation error, P. Błaszczyk (Poland), M. Katz (Israel), and D. Cherry (USA) arrive at the following conclusion: "Following Koetsier's hypothesis, it is reasonable to place it, rather, in the infinitesimal strand of the development of analysis, rather than the epsilontic strand.

After constructing the lower and upper sequences, Cauchy does write that the values of the latter "finiront par differer de ces premiers valeurs aussi peu que l'on voudra". That may sound a little bit epsilon/delta. Meanwhile, Leibniz uses language similar to Cauchy's: `Whenever it is said that a certain infinite series of numbers has a sum, I am of the opinion that all that is being said is that any finite series with the same rule has a sum, and that the error always diminishes as the series increases, so that it becomes as small as we would like ["ut fiat tam parvus quam velimus"].

Cauchy used epsilontics if and only if Leibniz did, over a century before him". [Błaszcyk, Katz, Sherry, p.18]

According to a Moscow researcher A.V. Dorofeeva writing about the mean value theorem, in Cauchy's works, "this conclusion is true only if the same $\delta$ can be picked out for all $x$, the fact whereof needs to be proved" [Dorofeeva, p. 48].

In 1985, a book of Bruno Belhoste entitled 'Cauchy. 1789 – 1857' [Belhoste, 1985] was published in Paris. Its translation [Belhoste,1997] into

Russian was published in 1997. Belhoste wrote in [Belhoste, p. 90] regarding Cauchy's proof of this Lagrange's theorem: "Instead of formula

$$f(x+i) - f(x) = pi + qi^2 + ri^3 + \ldots$$

that enabled Lacroix to present an increase of a function expandable to series and to define the differential, Cauchy proved the theorem on finite increments: if function $f$ is continuously differentiable between $x$ and $x+i$, then there exists such real positive number $\theta < 1$ that

$$f(x+i) - f(x) = i \cdot f'(x + \theta i).$$

He developed this formula from the inequation below using the theorem on intermediate values set forth in Analysis

$$\inf_{x \in [x_0, X]} f'(x) \leq \frac{f(X) - f(x_0)}{X - x_0} \leq \sup_{x \in [x_0, X]} f'(x). \qquad (*)$$

This inequation is correct for any continuous function (and hence, for differentiable function in the sense of Cauchy) between $x_0$ and $X$".

Please note that the theorem on intermediate values in The Course of Analysis of 1821 [Cauchy, 1821, p. 50] reads as follows: "*The Theorem on Continuous Function*. If function *f(x)*, a continuous function of variable *x* between limits *x=x₀, x=X* and *b,* is located between *f(x₀)* and *f(X),* then equation *f(x)=b* will always possess a solution located between *x₀* and *X.* "

Belhoste provides drawings to accompany Cauchy's theorems, much as we complement Lagrange's theorem with a function graph and feature a chord drawing extreme points together at lectures we read to students. However, you will find no drawing in the course of Cauchy; geometric interpretation of theorems is not mentioned anywhere either[14]. The statement provided by Belhoste is modern in its nature.

Thereafter, Belhoste continues: "The proof provided by Cauchy in 1823 only for functions continuously differentiable in $[x_0, X]$ made his new methods famous and made it possible to see the difference between a simple and uniform continuity.

However, his proof of the inequation (*) was based on a completely wrong assumption: if function $f$ is continuous (and therefore differentiable in the sense of Cauchy) between $x_0$ and $X$, and if

$\varepsilon$ is a positive number which is as small as we may wish, then, according to Cauchy, there exists such positive number $\delta$[15] that for any $i$ which is smaller than $\delta$ and for all $x$ between $x_0$ and $X$

$$f'(x) - \varepsilon \leq \frac{f(x+i) - f(x)}{i} \leq f'(x) + \varepsilon.$$

In fact, this inequation is true only for all $x$ located between $x_0$ and $X$, provided always that $f'$ is equicontinuous between these two numbers (or continuous in a closed bounded interval $[x_0, X]$). This error has proven that the lack of a clear distinction between continuity and uniform continuity was the weak point in the course of Cauchy. Nevertheless, the theorem on finite increments was consistently used and appeared to be the basic theorem in differential calculus" [Belhoste, 1985, p. 90–91].

It should be noted that it was the closed bounded interval that was meant by both Ampère and Cauchy. All examples illustrating this theorem were given for elementary functions that were uniformly continuous in a closed interval. Let us repeat Cauchy's words: "Function $f(x)$ is also said to always be a continuous

---

[14] There are no drawings in works of Cauchy, Lagrange, or Ampère. They will only appear in works of Lacroix [Lacroix, 1797], however, not to illustrate this theorem. Belhoste provides a modern geometric interpretation.

The author is thankful to S.S. Demidov for the following remark: *"Aged Lacroix certainly works in the manner of the 18th century. Therefore, he should not be regarded as going after Lagrange, Cauchy, and Ampère in terms of development of analysis. He merely did not master the manner introduced by Lagrange and followed by Cauchy and Ampère (not Lacroix): there should not be any drawings in the text – no reference to visualization!"* [Personal message] function of variable *x* in the vicinity of any particular value of such variable, if this function is continuous between two (even quite close) limits containing this given point" [Cauchy, 1823, p. 17].

[15] Note that Belhoste expressly provided that delta is chosen judging from epsilon, while Cauchy made no such express provision.

Cauchy did not use the language of $\varepsilon$–$\delta$[16] any more, not even in his late works. According to A.P. Yushkevich, "Cauchy's definition of continuity is as far from 'epsilontics' as his definition of limit" [Yushkevich, 1986, p. 69]. For the method to work, $\varepsilon$ and $\delta$ must be interrelated and have the structure of an interval (domain). In 1823, the understanding of a continuum was not yet developed enough for this. Let us also mention the standpoint of H. Putnam, which was as follows: "If the epsilon-delta methods had not been discovered, then infinitesimals would have been postulated entities (just as 'imaginary' numbers were for a long time). Indeed, this approach to the calculus enlarging the real number system–is just as consistent as the standard approach, as we know today from the work of Abraham Robinson. If the calculus had not been 'justified' Weierstrass style, it would have been 'justified' anyway" [Putnam, 1974].

Development of epsilontics was associated with the development of the concept of continuity. The remarkable similarity of Cauchy's and Bolzano's ideas had lead an English historian of mathematics, Ivor Grattan-Guinness, to a disputable idea of assimilation [Grattan-Guinness, 1970].

There are many examples in the history of science where the same ideas occurred to different scientists contemporaneously. One can disagree with Grattan-Guinness that such contemporaneous ideas were rather borrowed. The custom of previous problem statement could be so strong that it caused similar response of mathematicians working in different countries. This was the case with non-Euclidean geometry. This was the case with the concept of continuous function where Bolzano and Cauchy were based on Lagrange's ideas. This was the case with the concept of a surd number and continuity of continuum when Meray, Heine and Cantor contemporaneously offered similar concepts based on Cauchy's converging sequence criterion.

In 1868, 1869 and 1872, works of Charles Meray where he develops the theory of surd numbers with the help of a limit were published. The most complete statement of his theory is in the volume of 1872 [Méray, 1872]; comments can be found in works of Pierre Dugac [Dugac, 1973], [Dugac, 1972], [Sinkevich, 2012 c].

---

[16] The author of the article is responsible for this statement. All works of Cauchy are available.

During the XIX century, the necessity to express the relationship between epsilon and delta manifested itself more and more, but its functionality has formalized gradually. For example, Riemann, identifying improper integral

$$f_1(x + \delta) - f_1(x) = \int\limits_{x}^{x+\delta} f(x)dx = \delta f(x + \varepsilon\delta)$$

in the ε-neighborhood of the gap, said δ <ε. Hankel when defined limit used the estimate δ <ε [Hankel, p.195]. Eduard Heine, Ulisse Dini also enjoyed estimate δ <ε [Sinkevich 2012a, 2012b].

In 1854, Karl Weierstrass starts reading lectures in Berlin University. It was he who introduced such symbols as $\lim\limits_{n=\infty} p_n = \infty$ (published in 1856) [Yushkevich, 1986, p.76]).

Unfortunately, Weierstrass himself had never published or edited his lectures. In most cases, they came down to us in notes of his students. Eduard Heine bewailed in this regard: "Principles of Mr. Weierstrass are set forth directly in his lectures and in indirect spoken messages, and in handwritten copies of his lectures, and are pretty widely spread; however, their author's editions have never been published under the author's supervision, which damages their perceptual unity" [Heine, p.172]. But the main conception of ε–δ method formed in his Berlin lectures. According to A.P. Yushkevich, "Modern statement of differential calculus with its ε, δ-technique, wordings, and proofs, is reported to date back to lectures Weierstrass read in Berlin University, interpretations whereof were published by his students" [Yushkevich, 1977, p. 192].

The earliest known Weierstrass' text where the ε–δ technique is mentioned are differential calculus lecture notes made at a lecture read in the summer term of 1861 in Königlichen Gewerbeinstitut of Berlin. "The lecture notes were made by Weierstrass' student, H.A. Schwarz, and are now kept in Mittag-Leffler Institute in Sweden. Schwarz was 18 then, and he wrote these notes solely for himself, not to be published" [Yushkevich, 1977, p. 192]. Schwarz' notes were found and published by P. Dugac [Dugac, 1972]. It is in these notes that the definition of continuous function in the language of epsilontics appears for the first time: "If f (x) is function x and x is a defined value, then, on conversion of x into x+h, the function will change and will be f (x+h); difference f (x+h) − f (x) is used to be called the change received by the function by virtue of the fact that the argument convers from $x$ to $x + h$. If it is possible to determine such boundary $\delta$ for $h$ that for *all* values of $h$, the absolute value whereof is still smaller than $\delta$, $f(x+h) - f(x)$ becomes smaller than any arbitrarily small value of $\varepsilon$, then infinitesimal changes of function are said to correspond to infinitesimal changes of argument. Because a value is said to be able to become infinitely small, if its absolute value can become smaller than any arbitrary small value. If any function is such that infinitesimal changes of function correspond to infinitesimal changes of argument, then it is said to be a *continuous function* of argument or that it continuously changes along with its argument" [Yushkevich, 1977, p.189].

In 1872, Eduard Heine in «Die Elemente der Functionenlehre» [Heine, 1872] gave a definition of the limit function using Cantor's fundamental sequences. Every convergent sequence was represented as the sum of its limit and the elementary (decreasing) sequence [Heine, 1872, p.178]. On this basis, Heine formulates the condition of continuity [ibid, 182-183], the definition of uniform continuity in terms of $\varepsilon$-$\delta$, the theorem of uniformly continuous functions and as a method of proof of there was cover lemma. The theorem of uniformly continuity was necessary for intervals between irrational number and its limit, the rational number[17]. In the same year in Cantor's "Ueber die Ausdehnung eines Satzes der Theorie der trigonometrischen Reihen" there is a notion of a limit point. It was very useful and became widely known to mathematicians in Germany and Italy thanks to Schwarz. Ulysse Dini first used the concept of a limit point in its course "Fondamenti per la teoria delle funzioni di variabili reali". Theory of real number Dini expounded according to Dedekind. Defining the limit of the function for the case of the final argument, Dini writes: "For values $y$ to have a certain limit on the right and left of a finite number $a$, for example, on the right, it is necessary and sufficient that for any arbitrarily small positive number $\sigma$ there exists a positive number $\varepsilon$, such that the difference $y_{a+\varepsilon} - y_{a+\delta}$ between the value of $y$ at the point $x = a + \varepsilon$ i. e. $y_{a+\varepsilon}$, or any other value $y_{a+\delta}$ corresponding to the value of $x$ in $a + \delta$, values $y$ of $x$ between $a$ and $a + \varepsilon$ ($a$ excluded) was numerically smaller than $\sigma$ » [Dini, 1878, p 26.]. Dini had written that he was the first who give such definition [ibid]. Dini identified unilateral and bilateral limits to provide the definition of continuity and classification of discontinuities. When define uniformly continuous functions, Dini writes: "Let us now turn to the special functions that are continuous in the finite interval ($\alpha$, $\beta$), and above all, consider such function that for arbitrarily small non-zero positive number $\sigma$, for each $x$, possesses the particular value between $\alpha$ and $\beta$ ($\alpha$ and $\beta$ included) will be a special (particular) number ε, different from zero and positive, such that for all values of δ, numerically smaller ε, for which the point remains in the range, will be carried out in absolute value $f(a+\delta)-f(a)<\sigma$. But at the same value σ (which corresponds to the number ε) it is possible that the number ε, which is suitable for the number $a$, is not suitable for other points in the same range, so if it is necessary to decrease it; in addition, the question arises, as it happens when you are infinitely close to the point of discontinuity of the function, which is continuous only in a general sense, as well as continuous in the interval, as $x$ approaches special points, ε can be reduced to the limit, but it will never reach zero (which would be the lower limit value of ε). In other words, it is doubtful that in some cases the number ε, different from zero, may be used for all values of $x$ from α to β (α and β included), and therefore appropriate to distinguish different types of continuity in the range (α, β), namely, uniform continuity and nonuniform continuity; when for arbitrarily small positive number σ there exists a non-zero and positive number ε, such that for all values of δ, is numerically smaller than ε, at which point $x+\delta$ in the range of (α, β) (α and β included) will be performed in absolute value $f(x+\delta)-f(x)<\sigma$. As Cantor showed, if $f(x)$ is continuous in the interval from α to β, so for any number σ it always exists a number ε, which is the same at all points of interval. Then all of the foregoing is unnecessarily"[Dini, 1878, p. 46].

---

[17] It was formulated by Cantor, as Heine wrote. Now it is Cantor-Heine theorem.

In 1881, Austrian mathematician Otto Stolz (1842-1905) wrote an article "The importance of B. Bolzano in the history of calculus" [Stolz, 1881]. He argues for the importance of Bolzano's ideas in the development of infinitesimal analysis, based on the history of mathematical ideas of Lagrange, Cauchy, Duhamel, du Bois-Reymond, Weierstrass, Cantor, Schwarz and Dini. As he writes in the beginning, "Cauchy relied on infinitesimal calculus, abandoning the limits of the method of Lagrange, believing that only infinitesimal methods provide the necessary rigor. Clarity and elegance of its presentation facilitated widespread and the universal adoption of his course. Even found significant shortcomings, as time has shown, can be eliminated by the adoption of consistent principles based on the Cauchy arithmetic considerations. A few years before Cauchy these same views, sometimes substantially more fully developed by Bernard Bolzano (...). Hankel recognizes its priority over Cauchy in a proper understanding of the theory of infinite series. These ideas were continued Schwartz, Dini and Weierstrass" [Stolz, 1881, p. 255-256].

In 1885 O. Stolz published a textbook "Lectures on general arithmetic according to a new point of view" [Stolz, 1885], which sets out Weierstrass' analysis as a continuation of Cauchy's principles, in the «ε-δ» language.

In 1886 Weierstrass lectured on the theory of functions and used a notion of limit point when defined continuum [Weierstrass, p.72] and ε-δ when defined continuous and uniformly continuous functions [ibid, p. 73-74].

The legend that it was Cauchy who created the language of epsilontics appeared thanks to H. Lebesgue who wrote in his 'Lectures in integrating and searching for primal functions' of 1904: "For Cauchy, function $f(x)$ is continuous for $x_0$ if, regardless of the value of positive number ε, one can find such number $\eta(\varepsilon)$ that inequation $|h| \le \eta(\varepsilon)$ results in $|f(x_0 + h) - f(x_0)| \le \varepsilon$; function $f(x)$ is continuous in $(a, b)$ if correlation between ε and $\eta(\varepsilon)$ can be chosen regardless of $x_0$ for any $x_0$ in $(a, b)$" [Lebesgue, 1904, p.13]. In this connection, A.P. Yushkevich wrote: "In his famous work on integration theory, for some reason A. Lebesgue ascribes the definition of continuity of functions in the point stated in terms of epsilontics of early 20[th] century to Cauchy and describes this definition as classic. This is one of the numerous examples of modernization of assertions of authors of earlier days even by such outstanding mathematicians as H. Lebesgue was." [Yushkevich, 1986, p. 69]

Unfortunately, most historical errors were caused by the fact that authors did not turn to source materials. Instead, they believed a loose paraphrase of a third party who normally used modern language. We saw Belhoste's interpretation through supremum and infinum above, we saw as he added a geometric image, we saw interpretations by Lebesgue, Stolz, and others. In 1978, a reference book [Alexandrova, 1978] was published where article 'Limit' reads as follows: "The definition of a limit through $\varepsilon$ and $\delta$ was provided by Bolzano (1817) and thereafter, by Cauchy (1820)" [ibid., p.13]. As you and I have seen for ourselves, that is not so. Bolzano in 1817 and Cauchy in 1821 provided qualitative definitions of a limit and definitions of a continuous function in terms of increments; Cauchy used $\varepsilon$ and $\delta$ once when he improved Ampère's proof; however, Cauchy used $\varepsilon$ and $\delta$ as final assessments of and error, where $\delta$ did not depend on ε. Bolzano never used this technique. According to Weierstrass' lecture notes of 1861, it was Weierstrass who was the first to use the language of $\varepsilon$ and $\delta$ as a method.

In 1821, when Cauchy was writing his 'Course of Analysis', E. Heine was born in Berlin. Fifty-one years later, the latter stated the concept of uniform continuity. K. Weierstrass was 6 in 1821. It took about 40 years for him to start using epsilontics to the full extent.